\documentclass{gtmon_a}
\pdfoutput=1

\proceedingstitle{Groups, homotopy and configuration spaces (Tokyo
  2005)}
\conferencestart{5 July 2005}
\conferenceend{11 July 2005}
\conferencename{Groups, homotopy and configuration spaces, 
                in honour of Fred Cohen's 60th birthday}
\conferencelocation{University of Tokyo, Japan}

\editor{Norio Iwase}
\givenname{Norio}
\surname{Iwase}

\editor{Toshitake Kohno}
\givenname{Toshitake}
\surname{Kohno}

\editor{Ran Levi}
\givenname{Ran}
\surname{Levi}

\editor{Dai Tamaki}
\givenname{Dai}
\surname{Tamaki}

\editor{Jie Wu}
\givenname{Jie}
\surname{Wu}

\title[Twisted Alexander Polynomial]
{Twisted Alexander polynomials
and a\\partial order on the set of prime knots}

\author{Teruaki Kitano}
\givenname{Teruaki}
\surname{Kitano}
\address{
Department of Information Systems Science\\
Faculty of Engineering\\\newline
Soka University\\
Tokyo, 192-8577\\Japan
}
\email{kitano@soka.ac.jp}

\author{Masaaki Suzuki}
\givenname{Masaaki}
\surname{Suzuki}
\address{Department of Mathematics\\
Akita University\\\newline
Akita, 010-8502\\Japan
}
\email{macky@math.akita-u.ac.jp}


\volumenumber{13}
\issuenumber{}
\publicationyear{2008}
\papernumber{14}
\startpage{307}
\endpage{321}
\doi{}
\MR{}
\Zbl{}

\received{31 May 2006}
\revised{19 September 2007}
\accepted{26 September 2007}
\published{25 February 2008}
\publishedonline{25 February 2008}
\proposed{}
\seconded{}
\corresponding{}
\version{}

\dedicatory{Dedicated to Fred Cohen
on the occasion of his 60th birthday}

\makeatletter
\def\cnewtheorem#1[#2]#3{\newtheorem{#1}{#3}[section]
\expandafter\let\csname c@#1\endcsname\c@theorem}

\AtBeginDocument{\let\bar\wbar\let\tilde\wtilde\let\hat\what}
\let\Bbb\mathbb

\newtheorem{theorem}{Theorem}[section]
\cnewtheorem{proposition}[theorem]{Proposition}
\cnewtheorem{corollary}[theorem]{Corollary}

\theoremstyle{definition}
\cnewtheorem{definition}[theorem]{Definition}
\cnewtheorem{example}[theorem]{Example}

\theoremstyle{remark}
\cnewtheorem{remark}[theorem]{Remark}

\makeatother  
\makeautorefname{subsection}{Subsection}

\newcommand{\xbar}{\bar{x}}
\newcommand{\ybar}{\bar{y}}

\begin{document}

\begin{abstract}
We give a survey of some recent papers by the authors and Masaaki Wada
\cite{KS1,KS2,kswpaper} relating the twisted Alexander polynomial with
a partial order on the set of prime knots.  We also give examples and
pose open problems.
\end{abstract}

\begin{webabstract}
We give a survey of some recent papers by the authors and Masaaki Wada
relating the twisted Alexander polynomial with a partial order on the
set of prime knots.  We also give examples and pose open problems.
\end{webabstract}

\begin{asciiabstract}
We give a survey of some recent papers by the authors and Masaaki Wada
relating the twisted Alexander polynomial with a partial order on the
set of prime knots.  We also give examples and pose open problems.
\end{asciiabstract}

\maketitle

\section{Introduction}

 The twisted Alexander
polynomial was introduced by Lin \cite{Lin} and Wada \cite{Wada}
independently.  Lin defined it for a knot by using a regular Seifert
surface and Wada defined it for a finitely presentable group by using
Fox's free differential calculus.  More generally, Jiang and Wang
\cite{JW93-1} studied twisted topological invariants for a 3--manifold
by using representations of the fundamental group.  Following this
work, much research from many viewpoints has been carried out (or is
currently being carried out) related to this invariant.  See the
bibliography for more precise details.

In this paper, we give a survey of some recent papers by the authors
and Masaaki Wada \cite{KS1,KS2,kswpaper} relating the twisted
Alexander polynomial with a partial order on the set of prime knots.
We also give examples and pose open problems.

In \fullref{sect-def} we recall the definition of the twisted
Alexander polynomial for a finitely presentable group as given by Wada.
His definition is purely algebraic.

In \fullref{sect-surj} we state one of the fundamental results, which
gives a relation between the existence of a surjection of groups and
the twisted Alexander polynomial; see \cite{kswpaper} for more
details.  This result gives us a new criterion for the non-existence
of a surjective homomorphism between two groups.  The corresponding
fact about the (classical) Alexander polynomial is well known; see
Crowell and Fox \cite{CF}.

In \fullref{sect-order} we explain through examples how to use twisted
Alexander polynomials to determine the partial order on the set of
prime knots defined by the existence of a surjective homomorphism
between the knot groups.  Using the result of \fullref{sect-surj}, we
determined in \cite{KS1} which pairs of prime knots from Rolfsen's
knot table \cite{R} have this relation.  Rolfsen's knot table lists
all the prime knots of $10$ crossings or less.

In \fullref{sect-consider}, we consider some problems which arise
naturally from the results of \fullref{sect-order}.  The existence
problem of a surjective homomorphism between knot groups, or more
generally 3--manifold groups, is related with the theory of degree one
maps and the period of a knot.

In \fullref{sect-prob} we pose further problems for future study.

\section{Definition of the
twisted Alexander polynomial}\label{sect-def}

In this section, we briefly recall the algebraic definition of the
twisted Alexander polynomial, by using Fox's free differential
calculus.

Let $G$ be a finitely presentable group.
Choose a presentation of $G$:
\[
G =
\langle x_1,\ldots,x_u ~|~ r_1,\ldots,r_v \rangle .
\]
Suppose
$\alpha \co  G \to\langle\ t_1,\dots,t_l \ |\ t_i t_j=t_j
t_i\rangle\cong{\mathbb Z}^l$ is 
a surjective homomorphism to the free abelian group
with generators $t_1,\ldots,t_l$
and $\rho \co  G \to GL(n;R)$
a linear representation,
where $R$ is a unique factorization domain.
These maps naturally induce ring homomorphisms
$\tilde{\rho}$ and $\tilde{\alpha}$
from ${\mathbb Z}[G]$ to ${\Bbb Z}[GL(n;R)]=M(n;R)$ and
${\mathbb Z}[{t_1}^{\pm 1},\ldots,{t_l}^{\pm 1}]$ respectively,
where
$M(n;R)$ denotes the matrix algebra of degree $n$ over $R$.
Then
$\tilde{\rho}\otimes\tilde{\alpha}$
defines a ring homomorphism
\[
{\mathbb Z}[G]\to
M\left(n;R[{t_1}^{\pm 1},\ldots,{t_l}^{\pm 1}]\right).
\]
Let
$F_u$ be the free group on
generators $x_1,\ldots,x_u$ and let
$$
\Phi\co \Bbb Z[F_u]\to
M\left(n;R[{t_1}^{\pm 1},\ldots,{t_l}^{\pm 1}]\right)
$$
be the composite of the surjection
${\mathbb Z}[F_u]\to{\mathbb Z}[G]$
induced by the fixed presentation and the map
$\tilde{\rho}\otimes\tilde{\alpha}\co 
{\mathbb Z}[G]\to
M(n;R[{t_1}^{\pm 1},\ldots,{t_l}^{\pm 1}])$.

We define the $v \times u$ matrix $M$
whose $(i,j)$ component is the $n \times n$ matrix
$$
\Phi\left(\frac{\partial r_i}{\partial x_j}\right)
\in M\left(n;R[{t_1}^{\pm 1},\ldots,{t_l}^{\pm 1}]\right),
$$
where
${\partial}/{\partial x}$
denotes the Fox differential.
This matrix $M$ is called
the Alexander matrix of
the fixed presentation of $G$
associated to the representation $\rho$.

It is easy to see that
there is an integer $1 \leq j \leq u$
such that
$\det \Phi(x_j-1)\neq 0$.
For such $j$,
let us denote by $M_j$
the $v \times (u-1)$ matrix obtained from $M$
by removing the $j$-th column.
We regard $M_j$ as
an $n v \times n (u-1)$ matrix with coefficients in
$R[{t_1}^{\pm 1},\ldots,{t_l}^{\pm 1}]$.
Moreover, for an $n (u-1)$--tuple of indices
\[
I =
\left( i_1, i_2, \ldots, i_{n(u-1)} \right)
\quad
\left( 1 \leq i_1 < i_2 < \cdots < i_{n(u-1)} \leq n v \right), 
\]
we denote by $M_j^I$ the
$n (u-1) \times n (u-1)$ square matrix
consisting of the $i_k$-th rows of the matrix $M_j$,
where $k=1,2,\ldots,n(u-1)$.

Then the twisted Alexander polynomial of 
a finitely presented group $G$
for a representation $\rho\co  G \to GL(n;R)$
is defined to be a rational expression
$$
\Delta_{G,\rho}(t_1,\ldots,t_l)
=\frac{\gcd_I (\det M_j^I)}{\det \Phi(x_j-1)}.
$$
Here $\gcd_I (\det M_j^I)$ is the greatest common divisor
of $\{\det M_j^I\ |\ I\subset \{1,\dots,nv\}\}$.
It is shown that 
$\Delta_{G,\rho}(t_1,\ldots,t_l)$ is independent
of the choice of $j$ such that $\det \Phi(x_j-1) \neq 0$
up to a factor
$\epsilon {t_1}^{\varepsilon_1} \cdots {t_l}^{\varepsilon_l}$,
where $\epsilon \in R^{\times},\varepsilon_i \in {\mathbb Z}$.
Moreover,
we can check that Tietze transformations
on the presentation of $G$
do not affect the twisted Alexander polynomial.
Therefore the twisted Alexander polynomial
$\Delta_{G,\rho}(t_1,\ldots,t_l)$ is independent
of the choice of the presentation of $G$.
See Wada \cite{Wada} for details here.

\section{Twisted Alexander polynomials and
surjectivity of a group homomorphism}\label{sect-surj}

In this section,
we present the following theorem from \cite{kswpaper},
which is one of the fundamental theorems, 
for applications \cite{KS1,KS2} of twisted Alexander polynomials.

\begin{theorem}\label{mainsurjthm}
Let $G$ and $G'$ be finitely presentable groups and
$\alpha,\alpha'$ surjective homomorphisms
from $G,G'$ to ${\mathbb Z}^l$ respectively.
Suppose that there exists a surjective homomorphism
$\varphi \co  G \to G'$ such that
$\alpha = \alpha' \circ \varphi$.
Then $\Delta_{G,\rho}$
is divisible by $\Delta_{G',\rho'}$
for any representation $\rho' \co  G' \to GL(n ; R)$,
where $\rho = \rho' \circ \varphi$.
More precisely, 
the quotient of $\Delta_{G,\rho}$
by $\Delta_{G',\rho'}$ is a Laurent polynomial in
$R[t^{\pm1}_1,\dots,t^{\pm}_l]$.
\end{theorem}

The main motivation here is the following.  Let $G(K)$ be the knot
group $\pi_1(S^3 - K)$ of a knot $K$ in $S^3$.  For any knots $K,K'$,
if there exists a surjective homomorphism from $G(K)$ to $G(K')$, then
the Alexander polynomial of $K$ is divisible by that of $K'$.
Murasugi \cite{murasugilect} mentions that if there exists a
surjective homomorphism from a knot group $G(K)$ to the trefoil knot
group, then the twisted Alexander polynomial of $G(K)$ is divisible by
that of the trefoil knot group.  \fullref{mainsurjthm} is a
generalization of his result.

\fullref{mainsurjthm} is proved by two different methods in
\cite{KS1}.  The first one is purely algebraic using the definition of
the twisted Alexander polynomial and linear algebra.  The second one
uses Reidemeister torsion.  The twisted Alexander polynomial of a knot
group may be interpreted in terms of Reidemeister torsion in the same
way as the classical Alexander polynomial.  These invariants of knots
have been studied from the Reidemeister torsion viewpoint.  For
example, see Milnor \cite{Milnor}, Kitano \cite{Kitano} or Kirk and
Livingston \cite{KL}.

We now consider an easy algebraic situation as an example.

\begin{example}
Let $G = \langle 
x_1,\ldots, x_u | r_1,\ldots,r_v 
\rangle$ be a finitely presented group 
with the abelianization $\alpha \co  G \to {\mathbb Z}^l$ 
and $\rho \co  G \to GL(n;R)$ a representation of $G$. 
Another group $G'$ is defined by 
\[
G' = 
\langle 
x_1,\ldots, x_u | r_1,\ldots,r_v, s
\rangle 
\]
where $s$ is a word of $x_1,\ldots,x_u$. 
The natural projection $\pi \co  G \to G'$ is surjective. 
Suppose that the word $s$ belongs 
to $\ker \alpha$ and $\ker \rho$,
then there exists 
a surjective homomorphism $\alpha' \co  G' \to {\mathbb Z}^l$ 
and a representation $\rho' \co  G' \to GL(n;R)$ such that 
$\alpha = \alpha' \circ \pi$ and $\rho = \rho' \circ \pi$. 
By \fullref{mainsurjthm}, 
$\Delta_{G,\rho}$ is divisible by $\Delta_{G',\rho'}$. 
Here we will verify the divisibility by the definition. 
Let $M$ be the Alexander matrix of $G$ 
associated to $\rho$. 
Then the Alexander matrix $M'$ of $G'$ 
associated to $\rho'$ is obtained by adding $n$ rows  
\[
\left(
\Phi \left( \frac{\partial s}{\partial x_1} \right) ~
\Phi \left( \frac{\partial s}{\partial x_2} \right) ~
\cdots ~
\Phi \left( \frac{\partial s}{\partial x_u} \right)
\right) 
\] 
to $M$. 
Since the numerator of $\Delta_{G',\rho'}$ is 
the greatest common divisor of $\det {M'}_j^{I'}$, 
$\gcd_{I'} (\det {M'}_j^{I'})$ is a divisor of 
$\gcd_{I} (\det {M}_j^{I})$. 
Furthermore, 
the denominators of $\Delta_{G,\rho}$ and $\Delta_{G',\rho'}$ 
are the same. 
Hence $\Delta_{G,\rho}$ is divisible by $\Delta_{G',\rho'}$. 
\end{example}

\section{A partial order in the knot table}\label{sect-order}

In this section we consider prime knots.  Let $K$ be a prime knot and
$G(K)$ its knot group.  A partial order on the set of prime knots is
defined as follows.  For two prime knots $K_1,K_2$, we write $K_1 \geq
K_2$ if there exists a surjective group homomorphism from $G(K_1)$
onto $G(K_2)$.  It is known that the relation $\geq$ satisfies the
condition of a partial order on the set of prime knots.  The following
result \cite{KS1} can be obtained by applying Theorem 3.1 and using a
computer.

\begin{theorem}\label{mainorderthm}
The above partial order on the knots
in Rolfsen's table is given as below:
\[
\begin{array}{l}
8_5,8_{10},8_{15},8_{18},8_{19},8_{20},8_{21},
9_1,9_6,9_{16},9_{23},9_{24},9_{28},9_{40}, \\
10_5,10_9,10_{32},10_{40},10_{61},10_{62},10_{63},10_{64},
10_{65},10_{66},10_{76},10_{77}, 10_{78},\\
10_{82},10_{84},10_{85}, 10_{87},10_{98},10_{99},10_{103},10_{106},
10_{112},10_{114}, 10_{139},10_{140},\\
10_{141},10_{142},10_{143},10_{144}, 10_{159},10_{164}
\end{array}
\geq 3_1,
\]
\[
\begin{array}{l}
8_{18},9_{37},9_{40},
10_{58},10_{59},10_{60},10_{122},10_{136},
10_{137},10_{138}
\end{array}
\geq 4_1,
\]
\[
10_{74},10_{120},10_{122}
\geq 5_2 .
\]
\end{theorem}

First,
we explain how to show the non-existence of surjective homomorphisms
by using the twisted Alexander polynomials.
Then,
we construct explicitly a surjective homomorphism
for each pair of knots
which belongs to the above list.

\subsection{Non-existence of surjective
homomorphisms}\setobjecttype{subsec}\label{subsect-nonex}

We see how to apply
the twisted Alexander polynomials of knots
for the non-existence problem of surjective homomorphisms
between knot groups.

For many pairs of knots
which do not appear in \fullref{mainorderthm}, 
we can show easily that 
there exists no surjective homomorphism
between their knot groups
by using the classical Alexander polynomial. 
If the Alexander polynomial of $G(K)$
is not divisible by that of $G(K')$,
then there exists no surjective homomorphism
from $G(K)$ to $G(K')$.
However, for some cases,
it cannot be determined
by using only the Alexander polynomial 
whether or not there exists a surjective homomorphism
between knot groups.

\begin{example}
We explain how to find
whether or not there exists a surjective homomorphism
between $G(8_{11})$ and $G(3_1)$.
The classical Alexander polynomials of them are respectively
\begin{eqnarray*}
\Delta_{8_{11}} & = &
2 t^4 - 7 t^3 + 9 t^2 - 7 t + 2 ,\\
\Delta_{3_{1}} & = &
t^2 - t + 1.
\end{eqnarray*}
Because $\Delta_{8_{11}}$ does not divide $\Delta_{3_{1}}$ clearly,
there exists no surjective homomorphism
from $G(3_1)$ to $G(8_{11})$.
Therefore we can see 
\[
3_1 \not \geq 8_{11} .
\]
However, since $\Delta_{3_{1}}$ divides $\Delta_{8_{11}}$,
we cannot determine
whether or not there exists a surjective homomorphism
from $G(8_{11})$ to $G(3_{1})$.
\end{example}

In cases that we can not determine 
the non-existence of a surjection, 
we apply the twisted Alexander polynomial to the knot groups. 

We take a Wirtinger presentation of any knot group $G(K)$ 
as follows:
\[
G(K) =
\langle
x_1,x_2,\ldots, x_u \, | \, r_1, r_2,\ldots, r_{u-1}
\rangle .
\]
Here we note that generators are conjugate with each other. 
The abelianization 
\[
\alpha\co G(K) \to {\mathbb Z} \cong \langle t \rangle
\]
is obtained by mapping each generator $x_i$ to $t$. 
We consider the twisted Alexander polynomial
for a 2--dimensional unimodular representation over a finite field.
Now we fix a prime integer $p$
and take a representation
$$
\rho \co  G(K) \to SL(2;{\Bbb F}_p).
$$
Here ${\Bbb F}_p$ is the finite prime field ${\Bbb Z} / p {\Bbb Z}$.
Then 
we obtain the Alexander matrix $M$ associated to $\rho$
in $M((u-1)\times u;M(2;{\Bbb F}_p[t,t^{-1}]))$.
Further it always holds that $\det \Phi (x_1 - 1)\neq 0$
because 
\begin{align*}
\Phi(x_1-1)
&=\alpha(x_1)\otimes\rho(x_1)-E\\
&=t\rho(x_1)-E\\
\end{align*}
where $E$ is the identity matrix of degree 2. 
We write $M_1$ for the $(u-1) \times (u-1)$ matrix
obtained from $M$ by removing the first column.
Here $\Delta^{N}_{K,\rho}(t)$ denotes the determinant of $M_1$
and $\Delta^{D}_{K,\rho}(t)$ the determinant  of $\Phi (x_1 - 1)$.

In this case, the twisted Alexander polynomial of $G(K)$
for a representation
$\rho \co  G(K) \to SL(2;{\Bbb F}_p)$
is defined to be
\[
\Delta_{K,\rho} (t) =
\frac{\Delta^{N}_{K,\rho}(t)}{\Delta^{D}_{K,\rho}(t)}
=
\frac{\det M_1}{\det \Phi (x_1 - 1)}.
\]
Here we can prove the non-existence
of a surjective homomorphism between
the groups of any two knots
except for the pairs listed in \fullref{mainorderthm}.
A criterion for the non-existence is obtained
by applying \fullref{mainsurjthm} for knot groups.

\begin{corollary}\label{criterion}
Let $K_1$ and $K_2$ be two knots.
If there exists a representation\break
$\rho_2 \co  G(K_2) \to SL(2;{\Bbb F}_p)$
such that $\Delta^{N}_{K_1,\rho_1}(t)$
is not divisible by $\Delta^{N}_{K_2,\rho_2}(t)$
or\break
$\Delta^{D}_{K_1,\rho_1}(t)\neq \Delta^{D}_{K_2,\rho_2}(t)$
for any representation
$\rho_1 \co  G(K_1) \to SL(2;{\Bbb F}_p)$,
then there exists no surjective homomorphism from $G(K_1)$ onto $G(K_2)$.
\end{corollary}

By applying \fullref{criterion}
with the aid of a computer,
we can prove that there exists no surjective homomorphism
between the remaining pairs of knots.
All the twisted Alexander polynomials
which we use to check the non-existence of surjections
are listed in \cite{KS1}.

\begin{example}
We can show the non-existence of a surjective homomorphism
from $G(8_{11})$ to $G(3_1)$ 
by using the twisted Alexander polynomials of
$SL(2;{\Bbb F}_5)$--represent\-ations as follows.

First, we compute that
the numerators and denominators of
the twisted Alexander
polynomials of $G(8_{11})$
associated to all $SL(2;{\Bbb F}_5)$--representations.
We obtain the pairs of
$(\Delta^{N}_{8_{11},\rho_i}(t) ,
\Delta^{D}_{8_{11},\rho_i}(t) )$ $(i = 1,2,\ldots ,10)$ as follows: 
\begin{eqnarray*}
\lefteqn{(\Delta^{N}_{8_{11},\rho_i}(t) ,
\Delta^{D}_{8_{11},\rho_i}(t) ) =} \\
&& (2 t^6+2 t^5+t^4+4 t^3+t^2+2 t+2, t^2+t+1) ,\\
&& (3 t^6+2 t^5+4 t^4+4 t^3+4 t^2+2 t+3 ,t^2+4 t+1) ,\\
&& (t^8+t^6+t^2+1, t^2+1) ,\\
&& (2 t^8+2 t^7+4 t^6+t^5+4 t^4+t^3+4 t^2+2 t+2,
t^2+4 t+1) ,\\
&& (2 t^8+3 t^7+4 t^6+4 t^5+4 t^4+4 t^3+4 t^2+3 t+2,
t^2+t+1) ,\\
&& (4 t^8+3 t^6+t^4+3 t^2+4, t^2+1) ,\\
&& (4 t^8+t^7+t^6+3 t^4+t^2+t+4, t^2+4 t+1) ,\\
&& (4 t^8+2 t^7+t^5+2 t^4+t^3+2 t+4, t^2+3 t+1) ,\\
&& (4 t^8+3 t^7+4 t^5+2 t^4+4 t^3+3 t+4, t^2+2 t+1) ,\\
&& (4 t^8+4 t^7+t^6+3 t^4+t^2+4 t+4, t^2+t+1) .
\end{eqnarray*}
On the other hand,
for a certain $SL(2;{\Bbb F}_5)$--representation $\rho_0$,
the numerators and denominators of
the twisted Alexander
polynomials of $G(3_{1})$ is given by
\[
\Delta^{N}_{3_{1},\rho_0}(t) =
t^4 + 2 t^3 + 2t^2 + 2 t + 1, \quad
\Delta^{D}_{3_{1},\rho_0}(t) =
t^2 + 2 t + 1 .
\]
For any $i$, it is seen that
$\Delta^{N}_{8_{11},\rho_i}(t)$
is not divisible by $\Delta^{N}_{3_1,\rho_0}(t)$ or
$\Delta^{D}_{8_{11},\rho_i}(t)\neq \Delta^{D}_{3_1,\rho_0}(t)$.
Then there exists no surjective homomorphism from $G(8_{11})$ onto $G(3_1)$.
Therefore we obtain
\[
8_{11} \not \geq 3_1 .
\]
\end{example}

Using Alexander polynomials, we cannot determine the non-existence of
surjections between knot groups for $201$ pairs of knots.  However, we
can prove the non-existence of surjections by using the twisted
Alexander polynomials.  To prove it, we take $2,3,5,7,11,17$ as a
prime integer $p$ of $SL(2;{\Bbb F}_p)$--representation for
$26,50,81,33,10,1$ cases respectively.  All the data to check them are
shown in \cite{KS1}.

\subsection{Construction of surjective
homomorphisms}\setobjecttype{subsec}\label{subsect-const}

In this subsection, 
we explain how to construct a surjective homomorphism between
the groups of each pair of knots
which appears in the list of \fullref{mainorderthm}.

\begin{example}\label{eg-surj}
We show that there exist
surjective homomorphisms
$G(8_5) \to G(3_1)$ and $G(8_{18}) \to G(3_1)$.
The knot group of $8_5$, $8_{18}$ and $3_1$
admit Wirtinger presentations as follows:
\begin{eqnarray*}
G(8_5) &=&
\left\langle
\begin{array}{c}
y_1,y_2,y_3,y_4,\\
y_5,y_6,y_7,y_8
\end{array}
\right.
\, \left| \,
\begin{array}{c}
y_7 y_2 \ybar_7 \ybar_1 , \,
y_8 y_3 \ybar_8 \ybar_2 , \,
y_6 y_4 \ybar_6 \ybar_3 , \,
y_1 y_5 \ybar_1 \ybar_4 , \\
y_3 y_6 \ybar_3 \ybar_5 , \,
y_4 y_7 \ybar_4 \ybar_6 , \,
y_2 y_8 \ybar_2 \ybar_7
\end{array}
\right\rangle, \\
G(8_{18}) &=&
\left\langle
\begin{array}{c}
y_1,y_2,y_3,y_4,\\
y_5,y_6,y_7,y_8
\end{array}
\right.
\, \left| \,
\begin{array}{c}
y_4 y_1 \ybar_4 \ybar_2 , \,
y_5 y_3 \ybar_5 \ybar_2 , \,
y_6 y_3 \ybar_6 \ybar_4 , \,
y_7 y_5 \ybar_7 \ybar_4 , \\
y_8 y_5 \ybar_8 \ybar_6 , \,
y_1 y_7 \ybar_1 \ybar_6 , \,
y_5 y_8 \ybar_5 \ybar_7
\end{array}
\right\rangle, \\
G(3_1) &=&
\langle
x_1,x_2,x_3 \, | \,
x_3 x_1 \xbar_3 \xbar_2,\, x_1 x_2 \xbar_1 \xbar_3
\rangle ,
\end{eqnarray*}
where $\xbar = x^{-1},\ybar = y^{-1}$.
We define a map $\varphi$ from $G(8_5)$ to $G(3_1)$ as follows:
\[
\begin{array}{cccc}
\varphi(y_1) = x_3 , &
\varphi(y_2) = x_2 , &
\varphi(y_3) = x_1 , &
\varphi(y_4) = x_3 , \\
\varphi(y_5) = x_3 , &
\varphi(y_6) = x_2 , &
\varphi(y_7) = x_1 , &
\varphi(y_8) = x_3 .
\end{array}
\]
It is easy to check that this map $\varphi$ gives a homomorphism 
by computing the images of the relators. 
Moreover, it is clear that $\varphi$ is surjective 
by its definition. 
Then we obtain
\[
8_5 \geq 3_1 .
\]
Next,
we define another map $\varphi'$ from $G(8_{18})$ to $G(3_1)$ as follows:
\[
\begin{array}{cccc}
\varphi'(y_1) = x_1 , &
\varphi'(y_2) = x_2 , &
\varphi'(y_3) = x_1 , &
\varphi'(y_4) = x_3 , \\
\varphi'(y_5) = x_3 , &
\varphi'(y_6) = x_1 x_3 \xbar_1 , &
\varphi'(y_7) = x_3 , &
\varphi'(y_8) = x_1 .
\end{array}
\]
Similarly it can be seen that the map $\varphi'$ is also
a surjective group homomorphism.
Then we get 
\[
8_{18} \geq 3_1 .
\]
\end{example}

In \cite{KS1},
we constructed surjective group homomorphisms
for all pairs of knots in \fullref{mainorderthm}
explicitly.
This completes the proof of \fullref{mainorderthm}.

We remark that we could find
many surjective homomorphisms by using computer.
Once they are found by computer, 
it is easy to check by hand that 
they are surjective homomorphisms.

\section[Further results related to
\ref{mainorderthm}]{Further results related to
\fullref{mainorderthm}}\label{sect-consider}

In the previous section, we determined the partial order ``$\geq$'' in
Rolfsen's knot table.  These results lead to the following problems:
\begin{enumerate}
\item
Are three knots $3_1$, $4_1$ and $5_2$ minimal elements?
\item
Which relation can a pair of a periodic knot
and its quotient knot realize?
\item
Which relation can a degree one map between knot exteriors realize?
\end{enumerate}

In this section, we give partial answers to these problems.  First, we
prove that $3_1$ and $4_1$ are minimal elements.  To do this, we study
a surjective homomorphism from a fibered knot to another one.  Next,
we determine which relation a pair of a periodic knot and its quotient
knot realizes.  Finally, we study which relation given in
\fullref{mainorderthm} can be induced by a degree one map.  We do not
know the complete answer of the realizing problem for a degree one map.
Here we describe which surjection constructed in the previous section
is realized by a degree one map.

\subsection{Surjection between fibered knots}\setobjecttype{subsec}\label{fiberedknot}

In this subsection,
we prove the minimality of $3_1$ and $4_1$.

\begin{theorem}\label{fiberedthm}
$3_1$ and $4_1$ are minimal elements
under this partial ordering.
\end{theorem}

To prove this theorem,
we start to study a surjective homomorphism from 
a fibered knot group to another knot group.
Let $K_1$ and $K_2$ be knots in $S^3$.
We assume that $K_1$ is a fibered knot of genus $g_1$.
Here we obtain the following.

\begin{proposition}
If there exists a surjective homomorphism
$\varphi\co G(K_1) \rightarrow G(K_2)$,
then $K_2$ is also a fibered knot.
Further the genus $g_2$ of $K_2$ is
less than or equal to $g_1$.
\end{proposition}

\begin{proof}
By the result of Neuwirth \cite{Neuwirth} and
Stallings \cite{Stallings},
the commutator subgroup $[G(K_1),G(K_1)]$
of $G(K_1)$ is a free group of rank $2g_1$,
because $K_1$ is a fibered knot of genus $g_1$.
It is isomorphic to the fundamental group of its fiber surface.
Restricting this surjection $\varphi$ on $[G(K_1),G(K_1)]$,
we have a surjection
$$
\varphi|_{[G(K_1),G(K_1)]}\co [G(K_1),G(K_1)]
\rightarrow [G(K_2),G(K_2)].
$$
Here it is clear that $[G(K_2),G(K_2)]$ is
also a finitely generated group,
because it is the image of $[G(K_1),G(K_1)]$
by a surjection $\varphi$.
Hence the commutator subgroup $[G(K_2),G(K_2)]$ of $G(K_2)$
is also a free group of a finite rank.
By applying the result of \cite{Neuwirth}
and \cite{Stallings} again,
it is seen that $K_2$ is also a fibered knot.
The genus of $K_2$ is denoted by $g_2$.
Since $\varphi$ is a surjection between free groups,
it is clear that $g_2$ is less than or equal to $g_1$.
This completes the proof.
\end{proof}

As a corollary,
we can prove \fullref{fiberedthm} as follows.

\begin{proof}[Proof of \fullref{fiberedthm}]
We put $K=3_1$ or $4_1$. 
They are the fibered knots of genus 1.
Now we assume that $K'$ is a non-trivial prime knot
such that $K\geq K'$.
Here there exists a surjective homomorphism
$\varphi\co G(K)\rightarrow G(K')$ from the above assumption.
By the above proposition,
$K'$ is a fibered knot of genus 1, too.
Since any genus 1 fibered knot is $3_1$ or $4_1$,
$K'$ is $3_1$ or $4_1$.
It means that $3_1$ and $4_1$ are minimal elements.
\end{proof}

We remark that Silver and Whitten
studied the same result in \cite[Proposition 3.11]{silverwhitten}

\subsection{Period of a knot}

The periods of knots with up to $10$ crossings are listed in
\cite{KS2}.  To supplement \fullref{mainorderthm}, we have the
following theorem.

\begin{theorem}\label{periodthm}
The following relations are realized 
by pairs of a periodic knot and its quotient knot:
\[
\begin{array}{l}
8_5,8_{15},8_{19},8_{21},
9_1,9_{16},9_{28},9_{40}, \\
10_{61},10_{63},10_{64},
10_{66},10_{76},10_{78},10_{98},
10_{139},10_{141},
10_{142},10_{144}
\end{array}
\geq 3_1,
\]
\[
\begin{array}{l}
8_{18},
10_{58},10_{60},10_{122},10_{136},10_{138}
\end{array}
\geq 4_1,
\]
\[
10_{120}
\geq 5_2 .
\]
\end{theorem}

\subsection{Degree one maps}\setobjecttype{subsec}\label{subsect-degonemap}

Earlier we constructed surjections to prove \fullref{mainorderthm}.
In this subsection, we study which surjection is induced by a degree
one map.

First,
we recall the definition of a degree one map
in the case of a knot exterior.
Let $E(K_i)$ be oriented knot exteriors of $K_i$ 
in $S^3$ for $i=1,2$. 
A continuous map
\[
f\co (E(K_1), \partial E(K_1))\rightarrow (E(K_2),\partial E(K_2))
\]
is called a degree one map
if its induced map
\[
f_*\co H_3(E(K_1),\partial E(K_1);{\Bbb Z})
\rightarrow H_3(E(K_2),\partial E(K_2);{\Bbb Z})
\]
has a degree one.
It is known that
a degree one map
$f\co (E(K_1),\partial E(K_1))\to$\break $(E(K_2),\partial E(K_2))$
induces a surjective homomorphism $f_*\co G(K_1)\to G(K_2)$. 
See Hempel \cite{Hempel} for this fact. 

Let $f\co (E(K_1), \partial E(K_1))\rightarrow
(E(K_2),\partial E(K_2))$ be a degree one map.
By the definition of degree one map,
the induced map
\[
f_*\co H_3(E(K_1),\partial E(K_1);{\Bbb Z})\rightarrow H_3(E(K_2),
\partial E(K_2);{\Bbb Z})
\]
is an isomorphism. From the homology long exact sequence,
we have
$$
0\rightarrow H_3(E(K_i),\partial E(K_i);{\Bbb Z})
\overset{\partial}{\rightarrow} H_2(\partial E(K_i);{\Bbb Z})
\rightarrow H_2(E(K_i);{\Bbb Z})
\rightarrow \cdots
$$
Since each knot exterior $E(K_i)$ is homologically a circle,
then $H_2(E(K_i);{\Bbb Z})$ is vanishing.
Hence the above boundary map
$$
\partial\co  H_3(E(K_i),\partial E(K_i);{\Bbb Z})
\rightarrow H_2(\partial E(K_i);{\Bbb Z})
$$
is an isomorphism.
By the naturality, 
it is seen that 
\[
f\co (E(K_1),\partial E(K_1))\rightarrow
(E(K_2),\partial E(K_2))
\] 
is a degree one map if and only if
$f|_{\partial E(K_1)}\co 
\partial E(K_1)\rightarrow \partial E(K_2)$ is
a degree one map.
Since $\partial E(K_i)$ is a $2$--dimensional torus,
the degree of $f$ is determined by the determinant of
$f_*\co H_1(\partial E(K_1);{\Bbb Z})\cong {\Bbb Z}\oplus {\Bbb Z}
\rightarrow H_1(\partial E(K_2);{\Bbb Z})\cong{\Bbb Z}\oplus{\Bbb Z}$. 

Now we fix a basis of $H_1(\partial E(K_i);{\Bbb Z})$
so that
the first element is the meridian $m_i$ of $K_i$
and the second one is a longitude $l_i$ of $K_i$.
By such fixing bases of $H_1(\partial E(K_i);{\Bbb Z})$,
we can represent $f_*$ as a $2\times 2$--matrix.
The determinant of this matrix is the degree of 
$f\co \partial E(K_1)\rightarrow \partial E(K_2)$.
In our examples in \cite{KS1},
any $f_*$ maps the meridian of $K_1$ to that of $K_2$.
Then we have to compute only the image $f_*(l_1)$,
which can be written $f_*(l_1)=am_2+bl_2$.
Under this setting,
$f_*$ on the boundary can represented by the matrix
$\begin{pmatrix}
1 & a \\
0 & b\\
\end{pmatrix}$.
Hence, it is clear that
$f$ is a degree one map if and only if $b=\pm 1$.

By using the above argument,
we can check which surjective homomorphism
in \cite{KS1} is induced by a degree one map.

\begin{example}
Here we consider $G(8_5),G(8_{18})$ and $G(3_1)$ of 
\fullref{eg-surj}.
We check whether or not 
the surjective homomorphisms
which are constructed in \fullref{eg-surj}
are induced by a degree one map.

First, we check whether or not 
the map $\varphi \co  G(8_5) \to G(3_1)$
of \fullref{eg-surj} is induced
by a degree one map.
Pairs of
the meridians and the longitudes of $G(8_5)$ and $G(3_1)$
with respect to the presentations of \fullref{eg-surj}
are chosen as follows respectively:
\[
(y_1,
y_7 y_8 y_6 \ybar_1 y_3 y_4 y_2
\ybar_5 \ybar_1 \ybar_1 \ybar_1 \ybar_1 ),\quad
(x_3,
\xbar_2 \xbar_3 \xbar_1 x_3 x_3 x_3 ) .
\]
Clearly, the image of the meridian of $G(8_5)$ under $\varphi$
is the meridian of $G(3_1)$.
On the other hand, the image of the longitude of $G(8_5)$ under $\varphi$
is
\begin{eqnarray*}
\varphi ( \mbox{longitude of } 8_5) &=&
\varphi ( y_7 y_8 y_6 \ybar_1 y_3 y_4 y_2
\ybar_5 \ybar_1 \ybar_1 \ybar_1 \ybar_1 ) \\
&=& x_1 x_3 x_2 \xbar_3 x_1 x_3 x_2
\xbar_3 \xbar_3 \xbar_3 \xbar_3 \xbar_3 \\
&=& \xbar_3 \xbar_3 \xbar_3 x_1 x_3 x_2
\xbar_3 \xbar_3 \xbar_3 x_1 x_3 x_2 .\\
\end{eqnarray*}
Then
the image of the longitude of $8_5$ is equal 
to the $(-2)$ times of the longitude of $3_1$ 
in the first homology group.
Then this surjective homomorphism $\varphi$ is not
induced by a degree one map.

Next, we check whether or not 
the map $\varphi' \co  G(8_{18}) \to G(3_1)$
of \fullref{eg-surj} is induced
by a degree one map.
Similarly,
we fix pairs of
the meridians and the longitudes of $G(8_{18})$ and $G(3_1)$
with respect to the presentations of \fullref{eg-surj}:
\[
(y_1,
\ybar_4 y_5 \ybar_6 y_7 \ybar_8 y_1 \ybar_2 y_3
),\quad
(x_1,
\xbar_3 \xbar_1 \xbar_2 x_1 x_1 x_1 ) .
\]
Similarly, the image of the meridian of $G(8_{18})$
under $\varphi'$
is the meridian of $G(3_1)$.
Moreover, the image of the longitude of $G(8_{18})$
under $\varphi'$
is
\begin{eqnarray*}
\varphi' ( \mbox{longitude of } 8_{18})
&=&
\varphi' (\ybar_4 y_5 \ybar_6 y_7 \ybar_8 y_1 \ybar_2 y_3 ) \\
&=&
\xbar_3 x_3 x_1 \xbar_3 \xbar_1 x_3 \xbar_1 x_1 \xbar_2 x_1 \\
&=&
x_1 \xbar_3 \xbar_1 x_3 \xbar_2 x_1 . \\
\end{eqnarray*}
Then it is seen that
the image of the longitude of $8_{18}$ is equal 
to the longitude of $3_1$.
Therefore this surjective homomorphism $\varphi'$ is
induced by a degree one map.
\end{example}

We can check which surjection of \fullref{mainorderthm}
is induced by a degree one map.
Finally, we obtain the following as the result.

\begin{theorem}\label{degreeonethm}
The following relations are realized by degree one maps:
\[
8_{18},
10_5,10_9,10_{32},10_{40},10_{103},10_{106},
10_{112},10_{114},10_{159},10_{164}
\geq 3_1,
\]
\[
9_{37},9_{40} \geq 4_1,
\]
\[
10_{74},10_{122} \geq 5_2 .
\]
\end{theorem}

We do not know whether or not other relations 
are realized by degree one maps. 

\section{Problems}\label{sect-prob}

In this section, we present other problems connected with the partial
order $\geq$.

In \fullref{fiberedknot}, we proved that $3_1$ and $4_1$ are minimal.
However, it remains open whether or not $5_2$ is minimal.  So the
first problem is the following.

(1)\qua Characterize and determine the minimal knots 
under this partial order. 

As we see in \fullref{fiberedknot}, if there exists a surjective
homomorphism from the knot group of a fibered knot onto another knot
group, then its target knot is fibered.  We can restrict the partial
order to the set of prime fibered knot.

(2)\qua Determine the partial order on the set of fibered knots. 

The next problem is to study the relation between this partial
ordering and knot invariants.  There are some invariants to measure
and classify the complexity of knots.  The most fundamental invariant
is the crossing number.  So the following problem arises naturally.

(3)\qua If there exists a surjective homomorphism from $G(K_1)$ to
$G(K_2)$, then is the crossing number of $K_1$ greater than that of
$K_2$?

By \fullref{mainorderthm}, the answer is positive in case the crossing
number is smaller than or equal to $10$.  Moreover, we can check that
there exists no surjective homomorphism from the knot groups of knots
with up to $10$ crossings to that of alternating knots with $11$
crossings.  It still remains open for higher crossing cases.

The next problem also arises naturally.

(4)\qua If there exists a surjective homomorphism
from $G(K_1)$ to $G(K_2)$, 
then is the bridge number of $K_1$ 
greater than or equal to that of $K_2$?

The answer is also positive in case 
the crossing number is smaller than or equal to $10$.

(5)\qua If there exists a surjective homomorphism from 
$G(K_1)$ to $G(K_2)$ and both of $K_1$ and $K_2$ 
are hyperbolic knots, 
then is the volume  of $S^3-K_1$ 
greater than that of $S^3-K_2$?

We can calculate the hyperbolic volume for a given knot 
by SnapPea \cite{snappea}, 
or find its value on the web page KnotInfo \cite{knotinfo}. 
Thus we can make sure that 
every pair of knots in the list of 
\fullref{mainorderthm} satisfies this inequality. 
 
In \fullref{subsect-degonemap}, we considered which surjective
homomorphisms are induced by degree one map.  However, this was
restricted to the surjective homomorphisms which are constructed in
\cite{KS1}.  We do not know whether or not other pairs $K_1 \geq K_2$
can be realized by degree one maps.  There may exist another
surjective homomorphism other than we constructed, which is induced by
degree one map.  Therefore we should study the following problem.

(6)\qua Decide all pairs of knot in Rolfsen's table 
which admit degree one maps. 

We could not determine all geometric interpretations 
for the existences of surjective homomorphisms. 
However, Ohtsuki, Riley and Sakuma gave 
a systematic construction of surjective homomorphisms 
between 2--bridge link groups \cite{ORS}. 

\subsubsection*{Acknowledgements}

The authors would like to express their thanks to Professors Sadayoshi
Kojima and Dieter Kotschick for their useful comments and to Professor
Makoto Sakuma for explaining to us results on periods of knots.

\bibliographystyle{gtart}
\bibliography{link}

\end{document}